\documentclass[10pt,a4paper,reqno]{amsart}
\usepackage{amsfonts,amsthm,latexsym,amsmath,amssymb,amscd,amsmath, epsf}
\usepackage[utf8]{inputenc}
\usepackage{latexsym,amssymb,amsthm}
\usepackage{url,graphics} 
\usepackage[dvips]{graphicx,epsfig} 

\newtheorem{thm}{Theorem}
\newtheorem{proposition}[thm]{Proposition}
\newtheorem{lemma}{Lemma}

\newtheorem{example}{Example}
\newtheorem{remark}{Remark}

\newtheorem{conjecture}{Conjecture}

\def\q#1.{{\bf #1.}}
\renewcommand\geq{\geqslant}
\renewcommand\leq{\leqslant}
\newcommand{\cD}{\mathcal{D}}

\newcommand{\be}{\begin{equation}}
\newcommand{\ee}{\end{equation}}

\newcommand{\K}{\mathbb{K}}

\begin{document}
          \numberwithin{equation}{section}

          \title[A note on Fr\"{o}berg's conjecture]
          {A note on Fr\"{o}berg's conjecture for forms of equal degrees}

\author[G. Nenashev]{Gleb Nenashev}
\address{ Department of Mathematics,
   Stockholm University,
   S-10691, Stockholm, Sweden}
\email{nenashev@math.su.se}

\begin{abstract}

In this note by using  elementary  considerations, we settle Fr\"{o}berg's conjecture for a large number of cases, 
when all generators of ideals have the same degree.  
 \end{abstract}
\maketitle

 Let $S = \K[x_1,\ldots,x_n]$ be the polynomial ring in $n$ variables, where $\K$ is an algebraicly closed field of characteristic zero. 
The number of variables $n$ and the field $\K$ will be fixed throughout the whole paper.
 Denote by $S_d$  the $d$-th graded component of $S$, i.e., the linear space of all homogeneous  polynomials of degree~$d$ in $n$~variables.

 In~\cite{Fr}  R.\,Fr\"{o}berg formulated the following conjecture:
\begin{conjecture} 
\label{FrC}
Let $f_1,\ldots,f_z$ be generic forms of degrees $a_1,\ldots,a_z$ respectively. Set
$I = {<}f_1,\ldots,f_z{>}$. The Hilbert series of $S/I$ is given by:
$$HS_{S/I} (t) = \left[ \frac{\prod_{i=1}^z(1-t^{a_i})}  {(1-t)^{n}} \right],$$
where $[..]$ means that we truncate a real formal power series at its first negative term.
\end{conjecture}
He proved Conjecture~\ref{FrC} for $2$ variables and noticed that the left-hand side is bigger or equal than the right-hand side in the lexicographic sense. 
Later in~\cite{An}  D.\,J.\,Anick  proved Conjecture~\ref{FrC}  for $3$ variables. 
The conjecture is trivial, when $z\leq n$, and according to R.\,Stanley's result in~\cite{St} it is true for $z=n+1$. 
For $4$ variables, in~\cite{MM} J.\, Migliore and R.\,M.\,Miró-Roig proved  that any ideal generated by generic forms has weak Lefschetz property (strong L.P. is enough for a proof of the conjecture).
In this note we present some related  results in the case when all degrees $a_1,\ldots,a_z$ are the same.

Let $\cD_d$ be any nonempty class of forms of degree $d$ closed under the linear changes of coordinates. For example: $\cD_d=S_d$ or $\cD_d$ is the set of all $d$-th powers of linear forms.

We will work with the Hilbert function of an ideal; it is easy to convert it to Hilbert function of the quotient algebra, because the sum of dimensions of $m$-th~graded components of~$S/I$ and of $I$ is the dimension of $S_m$. 
 For $\cD_d$ and $z$, denote by $HF_{(\cD_d,z)}(m)$ the dimension of the $m$-th graded component  of an  ideal generated by $z$ generic forms from $\cD_d$; and denote by $HS_{(\cD_d,z)}(t)=\sum HF_{(\cD_d,z)}(m)t^m$ the Hilbert series of this ideal. 
In~\cite{HL} M.\,Hochster and D.\,Laksov found the values of $HF_{(S_d,z)}(d+1)$ for any $d$ and $z$. 
Below we generalize their result for the $(d+k)$-th graded component, but we miss $2\cdot dim(S_k)$ possible values of $z$.

\begin{thm}
\label{mainthm}
Let $d$ and $k$ be positive integers. Then
\begin{itemize}
\item for $z\leq \frac{dim(S_{d+k})}{dim(S_k)}-dim(S_k)$, \quad $HF_{(\cD_d,z)}(d+k)=z\cdot dim(S_k)$;
\item for $z\geq \frac{dim(S_{d+k})}{dim(S_k)}+dim(S_k)$, \quad $HF_{(\cD_d,z)}(d+k)= dim(S_{d+k})$.
\end{itemize}

\end{thm}

In~\cite{Au}  M.\,Aubry got the result of the first type, his result cover only thin set of cases ($d$ is  larger than some complicated function depends on $k$ and $z$). 
In~\cite{MM} J.\,Migliore and R.\,M.\,Miró-Roig also wrote similar result as a consequence of Anick's work, however their result only for small $z$ (the upper bound for $z$ depends only on $d,k$; doesn't depend on number of variables).

\smallskip

As a consequence of Theorem~\ref{mainthm} we get the following statement.
\begin{proposition}
\label{criteria}
Let $d$ and $z$ be positive integers.  If there exists $r$ such that
  $$\frac{dim(S_{d+r+1})}{dim(S_{r+1})}+dim(S_{r+1}) \leq z \leq \frac{dim(S_{d+r})}{dim(S_r)}-dim(S_r),$$
then the  Hilbert series of the ideal generated by $z$ generic forms from $\cD_d$ is given~by $$HS_{(\cD_d,z)}=\sum_{k=0}^\infty min(z\cdot dim(S_k),dim(S_{d+k}))t^{d+k}=\frac{1}{(1-t)^n}-\left[ \frac{(1-t^d)^z}{(1-t)^n}\right].$$
\emph{(}Recall that $\cD_d$ is any nonempty subset of $S_d$ closed under linear changes of coordinates.\emph{)}
\end{proposition}

Of course, all interesting cases correspond to $z\leq dim(S_d)$; otherwise  $HF_{(\cD_d,z)}(m)=dim(S_m)$ for $m\geq d$. 
Denote by  $p_d=\frac{\#\{{z\leq dim(S_d)}\text{ satisfying to Proposition 2}\}}{dim(S_d)}$ the "probability" that a given~${z\leq dim(S_d)}$  is covered by Proposition~\ref{criteria}.

\begin{example} For $n=5$ and $d=10$,
$dim(S_d)=1001$;   

$dim(S_1)=5$ and $\frac{dim(S_{d+1})}{dim(S_1)}=273$;

$dim(S_2)=15$ and $\frac{dim(S_{d+2})}{dim(S_2)}=\frac{364}{3}=121\frac{1}{3}$;

 $dim(S_3)=35$ and $\frac{dim(S_{d+3})}{dim(S_3)}=68$. 

\noindent Then the Hilbert series is given by Fr\"{o}berg's conjecture at least if number  of generators $z$ belongs to one of the following intervals:
\begin{itemize}
\item $z\geq 278$;
\item $268\geq z \geq 137$;
\item $106\geq z \geq 103$.
\end{itemize}
In other words, the Hilbert series is the standard one except possibly for $141=9+30+102$ cases. Thus  $$p_{10}=1-\frac{141}{1001}=0,859..$$
For larger $d$: $p_{15}=0,927..$; $p_{25}=0,968..$; $p_{40}=0,986..$
\end{example}

\begin{proposition}
\label{almostall}
For any fixed number of variables, the probability $p_d$ tends to $1$ as~$d\rightarrow +\infty$.

\end{proposition}

Proposition~\ref{almostall} means that Proposition~\ref{criteria} gives the criterium, which covers a huge number of nontrivial cases for large $d$. As a consequence, we get that Fr\"{o}berg's conjecture is true for many previously unknown  cases for large $d$ when the degrees of all forms are the same.

\smallskip
{\bf Acknowledgement.} I am grateful to Ralf~Fr\"{o}berg, Alessandro~Oneto for introducing me to this topic and to my supervisor Boris~Shapiro for his comments and help with this short text.

\section*{Proofs}
\begin{proof}[\bf Proof of Theorem~\ref{mainthm}] 

Fix $d$, $k$ and $\cD_d$. 
For a given $z$, define $a_z$ as the dimension of the $(d+k)$-th graded component of the intersection of two ideals generated by $z$ forms from $\cD_d$ and by one extra  form from $\cD_d$, which are generic.
 In other words, if $g_1,\ldots,g_z,g$ be generic forms, then
$$a_z:=dim\left( {<}g_1,\ldots,g_z{>}_{d+k}\cap {<}g{>}_{d+k}\right).$$
 \begin{lemma}
 If $a_{z+1}=a_z\neq 0$, then $a_z=dim(S_k)$ and  $HF_{(\cD_d,z)}(d+k)= dim(S_{d+k})$.
 \begin{proof}
 Consider generic forms $g_1,\ldots,g_z,g'_1,\ldots,g'_z$ and $g$ from $\cD_d$.
 
 We know that $$dim\left( {<}g_1,\ldots g_{z-1},g_z{>}_{d+k}\cap {<}g{>}_{d+k}\right)=a_z=$$
 $$= a_{z+1}=dim\left( {<}g_1,\ldots g_{z-1},g_z,g'_z{>}_{d+k}\cap {<}g{>}_{d+k}\right).$$
 We have $$dim\left( {<}g_1,\ldots g_{z-1},g_z{>}_{d+k}\cap {<}g{>}_{d+k}\right)=dim\left( {<}g_1,\ldots g_{z-1},g_z,g'_z{>}_{d+k}\cap {<}g{>}_{d+k}\right).$$
 The intersection in the left-hand side is a subspace of that in the right-hand side.
  Hence, they should coincide, we get
 $$ {<}g_1,\ldots g_{z-1},g_z{>}_{d+k}\cap {<}g{>}_{d+k}= {<}g_1,\ldots g_{z-1},g_z,g'_z{>}_{d+k}\cap {<}g{>}_{d+k}.$$
 Similarly, we have $$ {<}g_1,\ldots g_{z-1},g'_z{>}_{d+k}\cap {<}g{>}_{d+k}= {<}g_1,\ldots g_{z-1},g_z,g'_z{>}_{d+k}\cap {<}g{>}_{d+k} ,$$
 which implies $$ {<}g_1,\ldots g_{z-1},g_z{>}_{d+k}\cap {<}g{>}_{d+k}= {<}g_1,\ldots g_{z-1},g'_z{>}_{d+k}\cap {<}g{>}_{d+k}.$$
 Similarly, if we  change $g_{z-1}$  in right-hand side by the form $g'_{z-1}$, we get the same space. Repeating this procedure with $g_{z-2}$, $g_{z-3}$ and etc, we obtain 
    
 $$ {<}g_1,\ldots g_{z-1},g_z{>}_{d+k}\cap {<}g{>}_{d+k}= {<}g'_1,\ldots g'_{z-1},g'_z{>}_{d+k}\cap {<}g{>}_{d+k}.$$
 Hence for generic $g_1,\ldots,g_z,g\in\cD_d$, the linear space $V_g:={<}g_1,\ldots g_{z-1},g_z{>}_{d+k}\cap {<}g{>}_{d+k}$ depends only on $g$.
 
Fix any generic $g$, and choose a form $h\in V_g$. 
Hence for any generic $g_1,\ldots,g_z$, the form $h$ belongs to the ideal. 
For a linear coordinate transformation $A$, denote by $h_A$ the form $h$ after this  coordinate transformation.
Consider coordinate transformations $A_1,\ldots,A_b$ ($b$ is finite) such that the linear span of $h_{A_1},\ldots,h_{A_b}$ has the maximal dimension.

For generic $g_1,\ldots, g_z$ (generic with these $b$ coordinate transformations), the forms $h_{A_1},\ldots ,h_{A_b}$ belong to the ideal $I$ generated by $\{g_1,\ldots, g_z\}$.
 Hence, the linear span of $h_{A_1},\ldots,h_{A_b}$ belongs to the ideal $I$. Since this linear space has the maximal dimension, it is closed under the change of coordinates.

Hence, there is a nonempty linear space $H\subset S_{d+k}$ closed under the change of coordinates such that it  belongs to any ideal generated by generic $\{g_1,\ldots, g_z\}$. 
However $S_{d+k}$ has only one such subspace, namely $H=S_{d+k}$. Therefore, the $(d+k)$-th graded component of the ideal is the whole $S_{d+k}$. This proves the lemma. 

\end{proof}
 \end{lemma}

Let $z_0$ be the minimal $z$ such that $a_z\neq 0$, and $z_1$ be the minimal $z$ such that $a_{z_1}=dim(S_k)$.
By Lemma~1,  the dimension $a_z$ is strictly growing  between $z_0$ and~$z_1$, thus $$z_1-z_0\leq dim(S_k).$$
It is clear that \begin{itemize}
\item for $z\leq z_0$, the dimension $HF_{(\cD_d,z)}(d+k)=z\cdot dim(S_k)$;
\item for $z\geq z_1$, the dimension $HF_{(\cD_d,z)}(d+k)= dim(S_{d+k})$.
\end{itemize}
Since $z_0\leq \frac{dim(S_{d+k})}{dim(S_k)}$ and $z_1\geq \frac{dim(S_{d+k})}{dim(S_k)}$, we have $$z_1\leq z_0+dim(S_k) \leq \frac{dim(S_{d+k})}{dim(S_k)} +dim (S_k);$$
$$z_0\geq z_1-dim(S_k) \geq \frac{dim(S_{d+k})}{dim(S_k)} -dim (S_k),$$  
which gives the proof of the theorem.

\end{proof}

\begin{remark}
In fact, we proved  that $HF_{(\cD_d,z)}(d+k)=min(z\cdot dim(S_k), dim(S_{d+k}))$ except for at most possible $dim(S_k)$ cases for the values of~$z$. However, we don't know these $dim(S_k)$ values exactly.
\end{remark}
\bigskip

\begin{proof} [\bf Proof of Proposition~\ref{criteria}] 
By Theorem~\ref{mainthm}, we know  that $HF_{(\cD_d,z)}(d+r+1)=dim(S_{d+r+1})$ and $HF_{(\cD_d,z)}(d+r)=z\cdot dim(S_r)$.
 From the first claim, we get that the $(d+r+1)$-th graded component of the ideal is $S_{d+r+1}$, hence for $k\geq r+1$, the $(d+k)$-th graded component of ideal is $S_{d+k}$.

From the second claim, we get that for generic $g_1,\ldots,g_z$ from $\cD_d$, there are no $f_1,\ldots,f_z\in S_r$ (not all zeroes) such that $g_1f_1+\ldots+g_zf_z=0$. 
Hence, there are no such $f_1,\ldots,f_z\in S_k$, for $k\leq r$. 
Then for $k\leq r$, we have $HF_{(\cD_d,z)}(d+k)=z\cdot dim(S_k)$.
 Hence in this case, the whole Hilbert series is given by Fr\"{o}berg's conjecture.
 
\end{proof}

\begin{proof} [\bf Proof of Proposition~\ref{almostall}] 
Take an integer $k$. Then for large $d$, we know the Hilbert series for at leat
$$\sum_{r=0}^k \left(\left(\frac{dim(S_{d+r})}{dim(S_r)}-dim(S_r)\right)-\left(\frac{dim(S_{d+r+1})}{dim(S_{r+1})}+dim(S_{r+1})\right)\right)$$ 
different values of $z$. Then we have
$$1-p_d\leq dim(S_d)-\frac{\sum_{r=0}^k \left( \left(\frac{dim(S_{d+r})}{dim(S_r)}-dim(S_r)\right)-\left(\frac{dim(S_{d+r+1})}{dim(S_{r+1})}+dim(S_{r+1})\right)\right)}{dim(S_d)},$$ 

$$1-p_d\leq \frac{ \left( \frac{dim(S_{d+k+1})}{dim(S_{k+1})}\right)  }{dim(S_d)} +\frac{\sum_{r=0}^k \left( dim(S_r)+dim(S_{r+1})\right)}{dim(S_d)}.$$ 
The fist summand tends to $\frac{1}{dim(S_{k+1})}$ and the second one tends to zero as $d$ increases. Hence,  limsup of $(1-p_d)$ is at most  $\frac{1}{dim(S_{k+1})}$. 
Therefore,  $\lim_{d\to \infty}(1-p_d)=0$, because we have such a bound for any integer $k$.

\end{proof}

\end{document}